\documentclass[journal]{IEEEtran}
\usepackage{cite}
\usepackage[dvips]{graphicx}
\usepackage[cmex10]{amsmath}
\usepackage{amssymb}
\usepackage{arydshln}
\usepackage{multirow}
\usepackage{multicol}
\usepackage{color}

\newcommand{\R}{\mathbb{R}}
\newcommand{\N}{\mathbb{N}}
\newtheorem{theorem}{Theorem}
\newtheorem{lemma}{Lemma}
\newtheorem{remark}{Remark}

\newcommand{\change}[1]{#1}

\begin{document}

\title{Sampled-data implementation of derivative-dependent control using artificial delays}

\author{Anton~Selivanov and
        Emilia~Fridman,~\IEEEmembership{Senior Member,~IEEE}
\thanks{The authors are with the School of Electrical Engineering, Tel Aviv University, Israel (e-mail: antonselivanov@gmail.com; emilia@eng.tau.ac.il).}
\thanks{Supported by Israel Science Foundation (grant No. 1128/14). Preliminary results have been presented in \cite{Selivanov2017a}.}%
}

\markboth{}%
{}

\maketitle

\begin{abstract}
	We study a sampled-data implementation of linear controllers that depend on the output and its derivatives. First, we consider an LTI system of relative degree $r\ge 2$ that can be stabilized using $r-1$ output derivatives. Then, we consider PID control of a second order system. In both cases, the Euler approximation is used for the derivatives giving rise to a delayed sampled-data controller. Given a derivative-dependent controller that stabilizes the system, we show how to choose the parameters of the delayed sampled-data controller that preserves the stability under fast enough sampling. The maximum sampling period is obtained from LMIs that are derived using the Taylor's expansion of the delayed terms with the remainders compensated by appropriate Lyapunov-Krasovskii functionals. Finally, we introduce the event-triggering mechanism that may reduce the amount of sampled control signals used for stabilization. 
\end{abstract}

\section{Introduction}
Control laws that depend on output derivatives are used to stabilize systems with relative degrees greater than one. To estimate the derivatives, which can hardly be measured directly, one can use the Euler approximation $\dot y\approx (y(t)-y(t-\tau))/\tau$. This replaces the derivative-dependent control with the delay-dependent one \cite{Ilchmann2004,Niculescu2004,Karafyllis2008,Ramirez2017}. It has been shown in \cite{French2009} that \change{such} approximation preserves the stability if $\tau>0$ is small enough. Similarly, the output derivative in PID controller can be replaced by its Euler approximation. The resulting controller \change{was} studied in \cite{Okuyama2008} and \cite{Ramirez2016} using the frequency domain analysis. 

In this paper, we study \textit{sampled-data} implementation of the delay-dependent controllers. For double-integrators, this has been done in \cite{Liu2012} using complete Lyapunov-Krasovskii functionals with a Wirtinger-based term and in \cite{Seuret2015} via impulsive system representation and looped-functionals. Both methods lead to complicated linear matrix inequalities (LMIs) containing many decision variables. \change{In this paper, we obtain simpler LMIs for more general systems and prove their feasibility for small enough sampling periods.}

A simple Lyapunov-based method for delay-induced stabilization was proposed in \cite{Fridman2016,Fridman2017}. The key idea is to use the Taylor's expansion of the delayed terms with the remainders in the integral form that are compensated by appropriate terms in the Lyapunov-Krasovskii functional. This leads to simple LMIs feasible for small delays if the derivative-dependent controller stabilizes the system.

In this paper, we study sampled-data implementation of two types of derivative-dependent controllers. In Section~\ref{sec:2}, we consider an LTI system of relative degree $r\ge 2$ that can be stabilized using $r-1$ output derivatives. In Section~\ref{sec:3}, we consider PID control of a second order system. In both cases, the Euler approximation is used for the derivatives giving rise to a delayed sampled-data controller. Assuming that the derivative-dependent controller exponentially stabilizes the system with a decay rate $\alpha'>0$, we show how to choose the parameters of its sampled-data implementation that exponentially stabilizes the system with any decay rate $\alpha<\alpha'$ if the sampling period is small enough. The maximum sampling period is obtained from LMIs that are derived using the ideas of \cite{Fridman2016,Fridman2017}. Finally, we introduce the event-triggering mechanism that may reduce the amount of sampled control signals used for stabilization\change{ \cite{Astrom1999,Tabuada2007,Heemels2012,Heemels2013a,Selivanov2016c}}. In the preliminary paper \cite{Selivanov2017a}, we studied delayed sampled-data control for systems with relative degree two.

\textit{Notations:} $\N_0=\N\cup\{0\}$, $\mathbf{1}_r=[1,1,\ldots,1]^T\in\R^r$, $I_l\in\R^{l\times l}$ is the identity matrix, $\otimes$ stands for the Kronecker product\change{, $\lfloor x \rfloor=\max\{n\!\in\!\N\mid n\!\le\!x\}$ for $x\in\R$, $\operatorname{col}\{a_1,\ldots,a_r\}$ denotes the column vector composed from the vectors $a_1,\ldots,a_r$. For $p\in\R$, $f(h)=O(h^p)$ if there exist positive $M$ and $h_0$ such that $|f(h)|\le Mh^p$ for $h\in(0,h_0)$.}

\vspace{-.6em}
\subsection*{Auxiliary lemmas}
\begin{lemma}[Exponential Wirtinger inequality \cite{Selivanov2016f}]\label{lem:Wirt}
	Let $f\in\mathcal{H}^1(a,b)$ be such that $f(a)=0$ or $f(b)=0$. Then 
	\begin{multline*}\textstyle
	\int_a^be^{2\alpha t}f^T(t)Wf(t)\,dt\\ 
	\textstyle\le e^{2|\alpha|(b-a)}\frac{4(b-a)^2}{\pi^2}\int_a^be^{2\alpha t}\dot f^T(t)W\dot f(t)\,dt
	\end{multline*}
	for any $\alpha\in\R$ and $0\le W\in\R^{n\times n}$. 
\end{lemma}
\begin{lemma}[Jensen's inequality \cite{Solomon2013}]\label{lem:Jens}
	Let $\rho\colon[a,b]\to[0,\infty)$ and $f\colon[a,b]\to\R^n$ be such that the integration concerned is well-defined. Then for any $0<Q\in\R^{n\times n}$, 
	\begin{multline*}\textstyle
	\left[\int_a^b\rho(s)f(s)\,ds\right]^TQ\left[\int_a^b\rho(s)f(s)\,ds\right]\le\\ \textstyle
	\int_a^b\rho(s)\,ds\int_a^b\rho(s)f^T(s)Qf(s)\,ds.
	\end{multline*}
\end{lemma}
\section{Derivative-dependent control using discrete-time measurements}\label{sec:2}
Consider the LTI system 
\begin{equation}\label{LTI}
\begin{aligned}
\dot x(t)&=Ax(t)+Bu(t),\\
y(t)&=Cx(t),
\end{aligned}\qquad x\in\R^n, u\in\R^m, y\in\R^l
\end{equation}
with relative degree $r\ge2$, i.e., 
\begin{equation}\label{CB=0}
CA^iB=0,\quad i=0,1,\ldots,r-2,\quad CA^{r-1}B\neq 0. 
\end{equation}
Relative degree is how many times the output $y(t)$ needs to be differentiated before the input $u(t)$ appears explicitly. In particular, \eqref{CB=0} implies 
\begin{equation}\label{y'}
y^{(i)}=CA^ix,\quad i=0,1,\ldots, r-1. 
\end{equation}
To prove \eqref{y'}, note that it is trivial for $i=0$ and, if it has been proved for $i<r-1$, it holds for $i+1$: 
\begin{equation*}
\textstyle y^{(i+1)}=\left(y^{(i)}\right)'\stackrel{\eqref{y'}}{=}(CA^{i}x)'=CA^{i}[Ax+Bu]\stackrel{\eqref{CB=0}}{=}CA^{i+1}x. 
\end{equation*}

For LTI systems with relative degree $r$, it is common to look for a stabilizing controller of the form 
\begin{equation}\label{control:derivative}
u(t)=\bar{K}_0y(t)+\bar{K}_1\dot y(t)+\ldots+\bar{K}_{r-1}y^{(r-1)}(t) 
\end{equation}
\change{with $\bar{K}_i\in\R^{m\times l}$ for $i=0,\ldots,r-1$.}
\begin{remark}
	The control law \eqref{control:derivative} essentially reduces the system's relative degree from $r\ge2$ to $r=1$. Indeed, the transfer \change{matrix} of \eqref{LTI} has the form 
	\begin{equation*}\textstyle
	W(s)=\frac{\beta_rs^{n-r}+\cdots+\beta_n}{s^n+\alpha_1s^{n-1}+\cdots+\alpha_n}
	\end{equation*}
	with $\beta_r=CA^{r-1}B\neq0$. Taking $u(t)=K_0'u_0(t)+K_1'\dot u_0(t)+\cdots+K_{r-1}'u_0^{(r-1)}(t)$, one has 
	\begin{equation*}\textstyle
	\tilde y(s)=\frac{(\beta_{r}s^{n-r}+\cdots+\beta_n)(K_{r-1}'s^{r-1}+\cdots+K_0')}{s^n+\alpha_1s^{n-1}+\cdots+\alpha_n}\tilde u_0(s), 
	\end{equation*}
	\change{where $\tilde y$ and $\tilde u_0$ are the Laplace transforms of $y$ and $u_0$.} If $\beta_rK_{r-1}'\neq0$, the latter system has relative degree one. If it can be stabilized by $u_0=Ky$ then \eqref{LTI} can be stabilized by \eqref{control:derivative} with $\bar{K}_i=K_i'K$. 
\end{remark}

The controller \eqref{control:derivative} depends on the output derivatives, which are hard to measure directly. Instead, the derivatives can be approximated by the finite-differences 
\begin{equation*}
\begin{array}{l}
\dot y(t)\approx\frac{y(t)-y(t-\tau_1)}{\tau_1},\\
\ddot y(t)\approx\frac{1}{\tau_1}\left(\frac{y(t)-y(t-\tau_1)}{\tau_1}-\frac{y(t-\tau_1)-y(t-\tau_2)}{(\tau_2-\tau_1)}\right),\ldots\\
\end{array}
\end{equation*}
This leads to the delay-dependent control 
\begin{equation}\label{control:delay}
u(t)=K_0y(t)+K_1y(t-\tau_1)+\cdots+K_{r-1}y(t-\tau_{r-1}),
\end{equation}
where the gains $K_0,\ldots,K_{r-1}$ depend on the delays $0<\tau_1<\cdots<\tau_{r-1}$. If \eqref{LTI} can be stabilized by the derivative-dependent control \eqref{control:derivative}, then it can be stabilized by the delayed control \eqref{control:delay} with small enough delays \cite{French2009}. In this paper, we study the sampled-data implementation of \eqref{control:delay}: 
\begin{equation}\label{control:discrete}
\textstyle u(t)=K_0y(t_k)+\sum_{i=1}^{r-1} K_iy(t_{k}-q_ih),\quad t\in[t_k,t_{k+1}), 
\end{equation}
where $h>0$ is a sampling period, $t_k=kh$, $k\in\N_0$, are sampling instants, $0<q_1<\cdots<q_{r-1}$, $q_i\in\N$, are discrete-time delays, and $y(t)=0$ if $t<0$. 

In the next section, we prove that if \eqref{LTI} can be stabilized by the derivative-dependent controller \eqref{control:derivative}, then it can be stabilized by the delayed sampled-data controller \eqref{control:discrete} with a small enough sampling period $h$. Moreover, we show how to choose appropriate sampling period $h$, controller gains $K_0,\ldots,K_{r-1}$, and discrete-time delays $q_1,\ldots,q_{r-1}$. 
\subsection{Stability conditions}
Introduce the errors due to sampling 
\begin{equation*}
\begin{aligned}
\delta_0(t)&=y(t_k)-y(t),\\
\delta_i(t)&=y(t_k-q_ih)-y(t-q_ih),\\
\end{aligned}\quad t\in[t_k,t_{k+1}), k\in\N_0, 
\end{equation*}
where $i=1,\ldots,r-1$. Following \cite{Fridman2017}, we employ Taylor's expansion with the remainder in the integral form: 
\begin{equation*}\textstyle
y(t-q_ih)=\sum_{j=0}^{r-1}\frac{y^{(j)}(t)}{j!}(-q_ih)^j+\kappa_i(t), 
\end{equation*}
where 
\begin{equation*}\textstyle
\kappa_i(t)=\frac{(-1)^{r}}{(r-1)!}\int_{t-q_ih}^t(s-t+q_ih)^{r-1}y^{(r)}(s)\,ds. 
\end{equation*}
Combining these representations with \eqref{y'}, we rewrite \eqref{control:discrete} as 
\begin{equation}\label{control:errors}
u=[K_0, K]M\bar Cx+K_0\delta_0+K\delta+K\kappa, 
\end{equation}
where $\delta=\operatorname{col}\{\delta_1,\ldots,\delta_{r-1}\}$, $\kappa=\operatorname{col}\{\kappa_1,\ldots,\kappa_{r -1}\}$, 
\begin{equation}\label{MKC}
\begin{aligned}
M&=\left[\begin{smallmatrix}
I_l & 0 & 0 & \cdots & 0 \\
I_l & -q_1hI_l & \frac{(-q_1h)^2}{2!}I_l & \cdots & \frac{(-q_1h)^{r-1}}{(r-1)!}I_l \\
I_l & -q_2hI_l & \frac{(-q_2h)^2}{2!}I_l & \cdots & \frac{(-q_2h)^{r-1}}{(r-1)!}I_l \\
\vdots & \vdots & \vdots & \ddots & \vdots \\
I_l & -q_{r-1}hI_l & \frac{(-q_{r-1}h)^2}{2!}I_l & \cdots & \frac{(-q_{r-1}h)^{r-1}}{(r-1)!}I_l \\
\end{smallmatrix}\right], \\
K&=\begin{bmatrix}
K_1 & K_2 & \cdots & K_{r-1}
\end{bmatrix}, \quad \bar C=\left[\begin{smallmatrix}
C\\CA\\ \vdots \\CA^{r-1}
\end{smallmatrix}\right]. 
\end{aligned}
\end{equation}
The closed-loop system \eqref{LTI}, \eqref{control:discrete} takes the form 
\begin{equation}\label{ClosedLoop}
\begin{aligned}
\dot x&=Dx+BK_0\delta_0+BK\delta+BK\kappa, \\
D&=A+B[K_0, K]M\bar C.
\end{aligned}
\end{equation}
Using \eqref{y'}, the closed-loop system \eqref{LTI}, \eqref{control:derivative} can be written as
\begin{equation}\label{ClosedLoop:derivative}
\dot x=\bar{D}x,\quad \bar{D}=A+B[\bar{K}_0,\ldots,\bar{K}_{r-1}]\bar{C}. 
\end{equation}
Choosing 
\begin{equation}\label{gains}
[K_0,K_1,\ldots,K_{r-1}]=[\bar{K}_0,\ldots,\bar{K}_{r-1}]M^{-1}, 
\end{equation}
we obtain $D=\bar{D}$. (The Vandermonde-type matrix $M$ is invertible, since the delays $q_ih$ are different.) If \eqref{LTI}, \eqref{control:derivative} is stable, $\bar{D}$ must be Hurwitz and \eqref{ClosedLoop} will be stable for \change{zero} $\delta_0$, $\delta$, $\kappa$. The following theorem provides LMIs guaranteeing that $\delta_0$, $\delta$, and $\kappa$ \change{do not destroy} the stability of \eqref{ClosedLoop}. 
\begin{theorem}\label{th:1}
	Consider the LTI system \eqref{LTI} subject to \eqref{CB=0}.
	\begin{itemize}
		\item[(i)] For given sampling period $h>0$, discrete-time delays $0<q_1<\ldots<q_{r-1}$, controller gains $K_0, \ldots, K_{r-1}\in\R^{m\times l}$, and decay rate $\alpha>0$, let there exist positive-definite matrices $P\in\R^{n\times n}$, $W_0, W_i, R_i\in\R^{m\times m}$ ($i=1,\ldots,r-1$) such that\footnote{MATLAB codes for solving the LMIs are available at\\ \texttt{https://github.com/AntonSelivanov/TAC18}} $\Phi\le0$, where \change{$\Phi=\{\Phi_{ij}\}$ is the symmetric matrix} composed from
		\begin{equation*}\arraycolsep=1.4pt
		\begin{array}{rl}
		\Phi_{11}&=PD+D^TP+2\alpha P\\
		&\multicolumn{1}{r}{+h^2e^{2\alpha h}\sum_{i=0}^{r-1}(K_iCA)^TW_i(K_iCA),}\\
		\Phi_{12}&=\mathbf{1}_r^T\otimes PB,\quad\Phi_{13}=\mathbf{1}_{r-1}^T\otimes PB, \\
		\Phi_{14}&=(CA^{r-1}D)^TH,\quad \Phi_{24}=\mathbf{1}_r\otimes (CA^{r-1}B)^TH,\\
		\Phi_{22}&=-\frac{\pi^2}{4}\!\operatorname{diag}\{W_0,e^{-2\alpha q_1\!h}W_1,\ldots,e^{-2\alpha q_{r-1}\!h}W_{r-1}\!\},\\
		\Phi_{33}&=-(r!)^2\operatorname{diag}\Bigl\{\frac{e^{-2\alpha q_1h}}{(q_1h)^r}R_1,\ldots,\frac{e^{-2\alpha q_{r-1}h}}{(q_{r-1}h)^r}R_{r-1}\Bigr\},\\
		\Phi_{34}&=\mathbf{1}_{r-1}\otimes (CA^{r-1}B)^TH,\quad \Phi_{44}=-H\\
		\end{array}
		\end{equation*}		
		with $H=\sum_{i=1}^{r-1}(q_ih)^rK_i^TR_iK_i$ and $D$ defined in \eqref{ClosedLoop}. Then the delayed sampled-data controller \eqref{control:discrete} exponentially stabilizes the system \eqref{LTI} with the decay rate $\alpha$.
		\item[(ii)] Let there exist $\bar{K}_0,\ldots,\bar{K}_{r-1}\in\R^{m\times l}$ such that the derivative-dependent controller \eqref{control:derivative} exponentially stabilizes \eqref{LTI} with a decay rate $\alpha'$. Then, the delayed sampled-data controller \eqref{control:discrete} with $K_0,\ldots,K_{r-1}$ given by \eqref{gains} \change{with $M$ from \eqref{MKC}} and \change{$q_i=i\lfloor h^{\frac1r-1}\rfloor$} ($i=1,\ldots,r-1$)\footnote{\change{Note that $q_i=i\lfloor h^{\frac1r-1}\rfloor\in\N$ for small $h>0$}} exponentially stabilizes \eqref{LTI} with any given decay rate $\alpha<\alpha'$ if the sampling period $h>0$ is small enough. 
	\end{itemize}
\end{theorem}

{\em Proof} is given in Appendix~\ref{proof:th:1}. 

\begin{remark}\label{rem:h}
	Theorem~\ref{th:1}(ii) explicitly defines the controller parameters $K_0$, $K_i$, $q_i$ ($i=1,\ldots,r-1$), which depend on $h$. To find appropriate sampling period $h$, one should reduce $h$ until the LMIs from (i) start to be feasible. 
\end{remark}
\subsection{Event-triggered control}\label{sec:1:et}
\begin{figure}
  \centering
  \includegraphics[width=.69\linewidth]{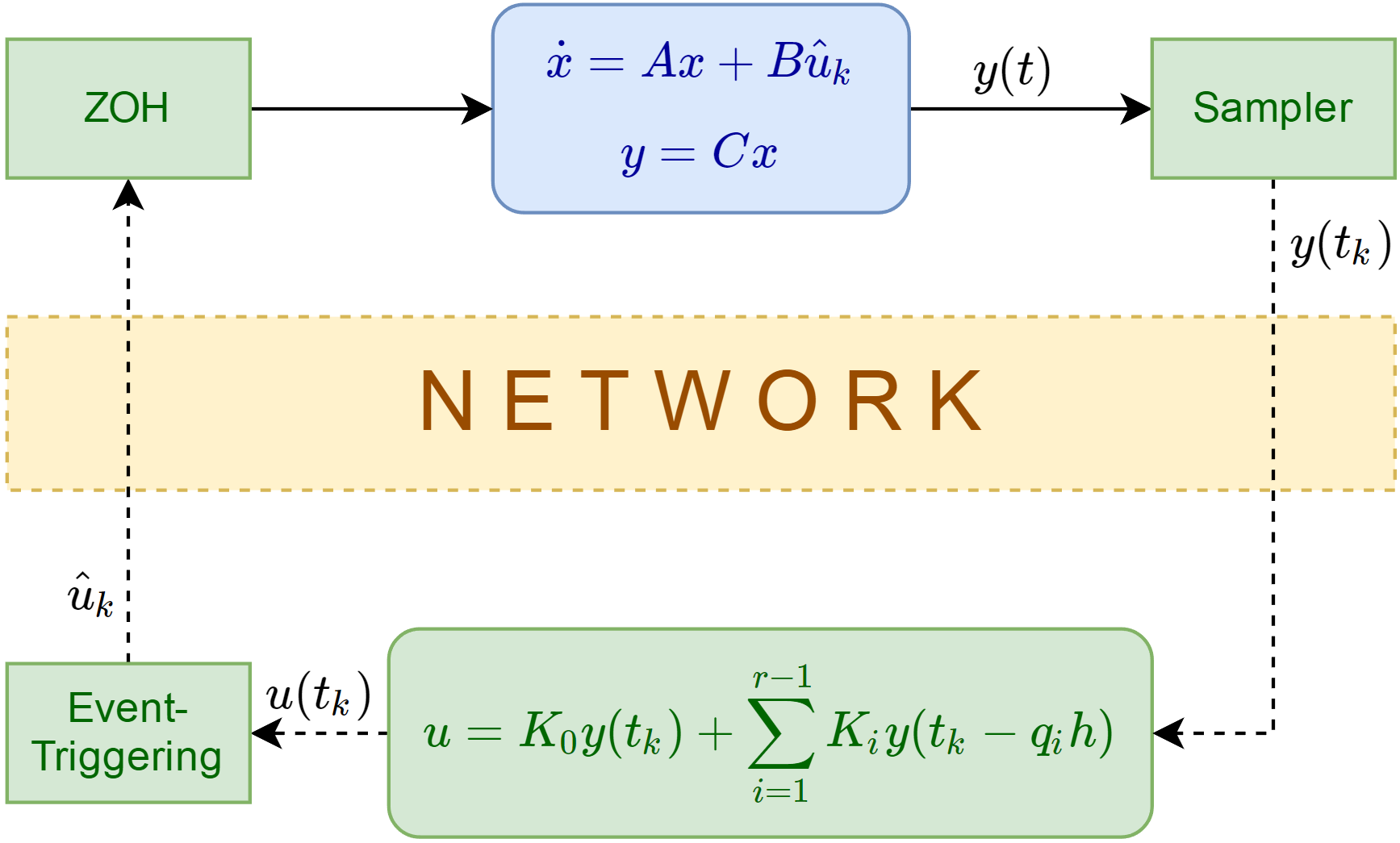}
  \caption{Event-triggering with respect to the control signal}
  \label{fig:ETC}  
\end{figure}

Event-triggered control allows to reduce the number of signals transmitted through a communication network \change{\cite{Astrom1999,Tabuada2007,Heemels2012,Heemels2013a,Selivanov2016c}}. The idea is to transmit the signal only when \change{it changes a lot}. The event-triggering mechanism for \textit{measurements} \change{was} implemented in \cite{Selivanov2017a} for the system \eqref{LTI}, \eqref{control:discrete} with relative degree $r=2$. Here, we consider the system with $r\ge2$ and introduce the event-triggering for \textit{control signals}, since the output event-triggering leads to complicated conditions (see Remark~\ref{rem:ET}). 

Consider the system (Fig.~\ref{fig:ETC})
\begin{equation}\label{LTI:ETC}
\begin{array}{ll}
\dot x(t)=Ax(t)+B\hat u_k,&t\in[t_k,t_{k+1}),\:k\in\N_0, \\
y(t)=Cx(t),&x\in\R^n, u\in\R^m, y\in\R^l, 
\end{array}
\end{equation}
where $\hat u_k=u(t_k)$ if $u(t_k)$ \change{from} \eqref{control:discrete} was transmitted and $\hat u_k=\hat u_{k-1}$ otherwise. The signal $u(t_k)$ is transmitted if its relative change since the last transmission is large enough, namely, if
\begin{equation}\label{ETC}
(u(t_k)-\hat u_{k-1})^T\Omega(u(t_k)-\hat u_{k-1})>\sigma u^T(t_k)\Omega u(t_k), 
\end{equation}
where $\sigma\in[0,1)$ and $0<\Omega\in\R^{m\times m}$ are event-triggering parameters. Thus, $\hat u_0=u(t_0)$ and 
\begin{equation}\label{ETControl}
\hat u_k=\left\{\begin{array}{ll}
u(t_k), & \text{\eqref{ETC} is true,}\\
\hat u_{k-1}, & \text{\eqref{ETC} is false.}
\end{array}\right.
\end{equation}

\begin{theorem}\label{th:2}
	Consider the system \eqref{LTI:ETC} subject to \eqref{CB=0}. For given sampling period $h>0$, discrete-time delays $0<q_1<\ldots<q_{r-1}$, controller gains $K_0, \ldots, K_{r-1}\in\R^{m\times l}$, event-triggering threshold $\sigma\in[0,1)$, and decay rate $\alpha>0$, let there exist positive-definite matrices $P\in\R^{n\times n}$, $\Omega, W_0, W_i, R_i\in\R^{m\times m}$ ($i=1,\ldots,r-1$) such that\footnote{MATLAB codes for solving the LMIs are available at\\ \texttt{https://github.com/AntonSelivanov/TAC18}} $\Phi_e\le0$, where 
	\begin{equation*}
	\begin{array}{l}\renewcommand*{\arraystretch}{.9}
	\Phi_e=\left[
	\begin{array}{cccc:cc}
	\multicolumn{4}{c:}{\multirow{4}{*}{$\Phi$}} & PB & \sigma([K_0, K]M\bar C)^T\Omega \\ 
	\multicolumn{4}{c:}{} & 0 & \sigma\mathbf{1}_r\otimes\Omega \\ 
	\multicolumn{4}{c:}{} & 0 & \sigma\mathbf{1}_{r-1}\otimes\Omega \\
	\multicolumn{4}{c:}{} & HCA^{r-1}\!B & 0 \\ \hdashline
	*&*&*&*& -\Omega & 0 \\
	*&*&*&*& 0 & -\sigma\Omega 
	\end{array} 
	\right]
	\end{array}
	\end{equation*}
	with $M$, $K$, $\bar C$ defined in \eqref{MKC} and $\Phi$, $H$ given in Theorem~\ref{th:1}. Then the event-triggered controller \eqref{control:discrete}, \eqref{ETC}, \eqref{ETControl} exponentially stabilizes the system \eqref{LTI:ETC} with the decay rate $\alpha$. 
\end{theorem}

{\em Proof} is given in Appendix~\ref{proof:th:2}. 
\begin{remark}\label{rem:ET}
The event-triggering mechanism \eqref{ETC}, \eqref{ETControl} is constructed with respect to the control signal. This allows to reduce the workload of a controller-to-actuator network. \change{To compensate the event-triggering error, we add \eqref{Sproc} to $\dot V$, which leads to two additional block-columns and block-rows in the LMI (confer $\Phi$ of Theorem~\ref{th:1} and $\Phi_e$ of Theorem~\ref{th:2}).} One can study the event-triggering mechanism with respect to the measurements by replacing $y(t_k)$, $y(t_k-q_ih)$ with $\hat y_k=y(t_k)+e_k$, $\hat y_{k-q_i}=y(t_k-q_ih)+e_{k-q_i}$ in \eqref{control:discrete}. This may reduce the workload of a sensor-to-controller network but \change{would require to add expressions similar to \eqref{Sproc} to $\dot V$ for each error $e_k, e_{k-q_1},\ldots, e_{k-q_{r-1}}$. This would lead to more complicated LMIs with two additional block-columns and block-rows for each error.} We study the event-triggering mechanism with respect to the control for simplicity. 
\end{remark}
\begin{remark}\label{rem:sigma}
	Taking $\Omega=\omega I$ with large $\omega>0$, one can show that $\Phi_e\le0$ and $\Phi\le0$ are equivalent for $\sigma=0$. This happens since the event-triggered control \eqref{control:discrete}, \eqref{ETC}, \eqref{ETControl} with $\sigma=0$ degenerates into periodic sampled-data control \eqref{control:discrete}. \change{Therefore, an appropriate $\sigma$ can be found by increasing its value from zero while preserving the feasibility of the LMIs from Theorem~\ref{th:2}.}
\end{remark}
\subsection{Example}\label{sec:tripint}
Consider the triple integrator $\dddot{y}=u$, which can be presented in the form \eqref{LTI} with 
\begin{equation}\label{integrator}
\renewcommand*{\arraystretch}{.8}
A\!=\!\left[\begin{array}{ccc}
0 & 1 & 0\\ 
0 & 0 & 1\\
0 & 0 & 0
\end{array}\right],\:
B=\left[
\begin{array}{c}
0\\ 
0\\
1
\end{array}
\right],\:
C=
\left[
\begin{array}{ccc}
1 & 0 & 0
\end{array}
\right]. 
\end{equation}
These parameters satisfy \eqref{CB=0} with $r=3$. The derivative-dependent control \eqref{control:derivative} with 
\begin{equation*}
\bar{K}_0=-2\times 10^{-4},\quad\bar{K}_1=-0.06,\quad\bar{K}_2=-0.342
\end{equation*}
stabilizes the system \eqref{LTI}, \eqref{integrator}. The LMIs of Theorem~\ref{th:1} are feasible for
\begin{equation*}
\begin{array}{c}
h=0.044,\quad q_1=30,\quad q_2=60,\quad \alpha=10^{-3}, \\
K_0\approx-0.265,\quad K_1\approx 0.483,\quad K_2\approx -0.219, 
\end{array}
\end{equation*}
where $K_i$ are calculated using \eqref{gains}. Therefore, the delayed sampled-data controller \eqref{control:discrete} also stabilizes the system \eqref{LTI}, \eqref{integrator}. 

Consider now the system \eqref{LTI:ETC}, \eqref{integrator}. The LMIs of Theorem~\ref{th:2} are feasible for $h=0.042$, $\sigma=2\times10^{-3}$ with the same control gains $\bar K_0$, $\bar K_1$, $\bar K_2$, delays $q_1$, $q_2$, and decay rate $\alpha$. Thus, the event-triggered control \eqref{control:discrete}, \eqref{ETC}, \eqref{ETControl} stabilizes the system \eqref{LTI:ETC}, \eqref{integrator}. Performing numerical simulations for $10$ randomly chosen initial conditions $\|x(0)\|_{\infty}\le1$, we find that the event-triggered control \eqref{control:discrete}, \eqref{ETC}, \eqref{ETControl} requires to transmit on average $455.6$ control signals during $100$ seconds. The amount of transmissions for the sampled-data control \eqref{control:discrete} is given by $\lfloor\frac{100}{h}\rfloor+1=2273$. Thus, the event-triggering mechanism reduces the workload of the controller-to-actuator network by almost $80\%$ preserving the decay rate $\alpha$. Note that $\sigma>0$ leads to a smaller sampling period $h$. Therefore, the event-triggering mechanism requires to transmit more measurements \change{through sensor-to-controller network}. However, the \change{total workload of both} networks is reduced by over $37\%$. 
\section{Event-triggered PID control}\label{sec:3}
Consider the scalar system 
\begin{equation}\label{PIDplant}
\ddot y(t)+a_1\dot y(t)+a_2y(t)=bu(t)
\end{equation}
and the PID controller
\begin{equation}\label{PID}
u(t)=\bar{k}_py(t)+\bar{k}_i\int_0^ty(s)\,ds+\bar{k}_d\dot y(t). 
\end{equation}
Here, we study sampled-data implementation of the PID controller \eqref{PID} that is obtained using the approximations
\begin{equation*}
\begin{array}{l}
\int_0^ty(s)\,ds\approx\int_0^{t_k}y(s)\,ds\approx h\sum_{j=0}^{k-1}y(t_j),\\
\dot y(t)\approx\dot y(t_k)\approx\frac{y(t_k)-y(t_{k-q})}{qh}, 
\end{array}\quad t\in[t_k,t_{k+1}), 
\end{equation*}
where $h>0$ is a sampling period, $t_k=kh$, $k\in\N_0$, are sampling instants, $q\in\N$ is a discrete-time delay, and $y(t_{k-q})=0$ for $k<q$. Substituting these approximations into \eqref{PID}, we obtain the sampled-data controller 
\begin{equation}\label{PID:SD}
\begin{array}{r}
u(t)=k_py(t_k)+k_ih\sum_{j=0}^{k-1}y(t_j)+k_dy(t_{k-q}),\\
t\in[t_k,t_{k+1}),\: k\in\N_0, 
\end{array}
\end{equation}
with $y(t_{k-q})=0$ for $k<q$ and 
\begin{equation}\label{PIDk}
\textstyle k_p=\bar{k}_p+\frac{\bar{k}_d}{qh},\quad k_i=\bar{k}_i,\quad k_d=-\frac{\bar{k}_d}{qh}. 
\end{equation}

Similarly to Section~\ref{sec:1:et}, we introduce the event-triggering mechanism to reduce the amount of transmitted control signals. Namely, we consider the system 
\begin{equation}\label{PIDplant1}
\ddot y(t)+a_1\dot y(t)+a_2y(t)=b\hat{u}_k,\quad t\in[t_k,t_{k+1}),\ k\in\N_0,  
\end{equation}
where $\hat u_k$ is the event-triggered control: $\hat u_0=u(t_0)$, 
\begin{equation}\label{ETCon}
\hat u_k=\left\{
\begin{array}{ll}
u(t_k),&\text{if \eqref{ETC:PID} is true},\\
\hat u_{k-1},&\text{if \eqref{ETC:PID} is false},
\end{array}
\right.
\end{equation}
with $u(t)$ from \eqref{PID:SD} and the event-triggering condition 
\begin{equation}\label{ETC:PID}
(u(t_k)-\hat u_{k-1})^2>\sigma u^2(t_k). 
\end{equation}
Here, $\sigma\in[0,1)$ is the event-triggering threshold. 
\begin{remark}
	We consider the event-triggering mechanism with respect to the control signal, since the event-triggering with respect to the measurements $\hat y_k=y(t_k)+e_k$ leads to an accumulating error in the integral term: 
	\begin{equation*}\textstyle
	\int_0^{t_k}y(s)\,ds\approx h\sum_{j=0}^{k-1}\hat y_j=h\sum_{j=0}^{k-1}y(t_j)+h\sum_{j=0}^{k-1}e_j. 
	\end{equation*}
\end{remark}
\subsection{Stability conditions}
To study the stability of \eqref{PIDplant1} under the event-triggered PID control \eqref{PID:SD}, \eqref{ETCon}, \eqref{ETC:PID}, we rewrite the closed-loop system in the state space. Let $x_1=y$, $x_2=\dot y$, and 
\begin{equation*}\textstyle
x_3(t)=(t-t_k)y(t_k)+h\sum_{j=0}^{k-1}y(t_j),\quad t\in[t_k,t_{k+1}). 
\end{equation*}
Introduce the errors due to sampling
\begin{equation*}
\begin{array}{l}
v(t)=x(t_k)-x(t)\\
\delta(t)=y(t_{k-q})-y(t-qh),
\end{array}\quad t\in[t_k,t_{k+1}),\ k\in\N_0. 
\end{equation*}
Using Taylor's expansion for $y(t-qh)$ with the remainder in the integral form, we have  
\begin{equation*}
y(t_{k-q})=y(t-qh)+\delta(t)\textstyle =y(t)-\dot y(t)qh+\kappa(t)+\delta(t), 
\end{equation*}
where 
\begin{equation*}\textstyle
\kappa(t)=\int_{t-qh}^t(s-t+qh)\ddot y(s)\,ds. 
\end{equation*}
Using these representations in \eqref{PID:SD}, we obtain 
\begin{equation}\label{uk}\arraycolsep=1.4pt
\begin{array}{rl}
u(t_k)&=k_px_1(t_k)+k_ix_3(t_k)+k_dy(t_{k-q})\\
&=[k_p+k_d, -qhk_d, k_i]x+[k_p, 0, k_i] v+k_d(\kappa+\delta). 
\end{array}
\end{equation}
Introduce the event-triggering error $e_k=\hat u_k-u(t_k)$. Then the system \eqref{PIDplant1} under the event-triggered PID control \eqref{ETCon}, \eqref{ETC:PID}, \eqref{uk} can be presented as 
\begin{equation}\label{plant}
\begin{array}{l}
\dot x=Ax+A_vv+Bk_d(\kappa+\delta)+Be_k,\\
y=Cx,
\end{array}
\end{equation}
for $t\in[t_k,t_{k+1})$, $k\in\N_0$, where 
\begin{equation}\label{PIDparameters}
\begin{array}{l}
A=\begin{bmatrix}
0 & 1 & 0 \\ -a_2+b(k_p+k_d) & -a_1-qhbk_d & bk_i \\ 1 & 0 & 0
\end{bmatrix},\\
A_v=\begin{bmatrix}
0 & 0 & 0 \\
bk_p & 0 & bk_i \\
1 & 0 & 0
\end{bmatrix}\!,\: 
B=\begin{bmatrix}
0 \\ b \\ 0
\end{bmatrix}\!,\:C=\begin{bmatrix}
1 & 0 & 0
\end{bmatrix}. 
\end{array}
\end{equation}
\change{Note that the ``integral'' term in \eqref{PID:SD} requires to introduce the error due to sampling $v$ that appears in \eqref{plant} but was absent in~\eqref{ClosedLoop}. The analysis of $v$ is the key difference between Theorem~\ref{th:2} and the next result.}
\begin{theorem}\label{th:3}
	Consider the system \eqref{PIDplant1}. 
	\begin{itemize}
		\item[(i)] For given sampling period $h>0$, discrete-time delay $q>0$, controller gains $k_p$, $k_i$, $k_d$, event-triggering threshold $\sigma\in[0,1)$, and 
		decay rate $\alpha>0$, let there exist positive-definite matrices $P, S\in\R^{3\times3}$ and nonnegative scalars $W$, $R$, $\omega$ such that\footnote{MATLAB codes for solving the LMIs are available at\\ \texttt{https://github.com/AntonSelivanov/TAC18}} $\Psi\le0$, where \change{$\Psi=\{\Psi_{ij}\}$ is the symmetric matrix} composed from
		\begin{equation*}
		\begin{array}{l}
		\Psi_{11}=PA+A^TP+2\alpha P+\left[\begin{smallmatrix}
		0 & 0 & 0 \\ 0 & 1 & 0 \\ 0 & 0 & 0
		\end{smallmatrix}\right]Wk_d^2h^2e^{2\alpha h},\\
		\Psi_{12}=PA_v\sqrt{h},\quad\Psi_{13}=\Psi_{14}=\Psi_{15}=PB,\\
		\Psi_{16}=\begin{bmatrix}
		k_p\!+\!k_d\\-qhk_d \\ k_i
		\end{bmatrix}\omega\sigma,\quad\Psi_{26}=\begin{bmatrix}
		k_p\\ 0 \\ k_i
		\end{bmatrix}\omega\sigma\sqrt{h},\\
		\Psi_{17}=A^TG,\quad\Psi_{22}=-\frac{\pi^2}{4}Sh,\quad\Psi_{27}=A_v^TG\sqrt{h},\\
		\Psi_{36}=\Psi_{46}=\omega\sigma,\quad\Psi_{37}=\Psi_{47}=\Psi_{57}=B^TG,\\
		\Psi_{33}=-W\frac{\pi^2}{4}e^{-2\alpha qh},\quad\Psi_{44}=-R\frac{4}{(qh)^2}e^{-2\alpha qh},\\
		\Psi_{55}=-\omega,\quad\Psi_{66}=-\omega\sigma,\quad\Psi_{77}=-G,\\
		G=h^2e^{2\alpha h}S+\left[\begin{smallmatrix}
		0 & 0 & 0 \\ 0 & 1 & 0 \\ 0 & 0 & 0
		\end{smallmatrix}\right]R k_d^2(qh)^2
		\end{array}
		\end{equation*}
		with $A$, $A_v$, $B$, $C$ given in \eqref{PIDparameters}. Then, the event-triggered PID controller \eqref{PID:SD}, \eqref{ETCon}, \eqref{ETC:PID} exponentially stabilizes the system \eqref{PIDplant1} with the decay rate $\alpha$. 
		\item[(ii)] Let there exist $\bar{k}_p$, $\bar{k}_i$, $\bar{k}_d$ such that the PID controller \eqref{PID} exponentially stabilizes the system \eqref{PIDplant} with a decay rate $\alpha'$. Then, the event-triggered PID controller \eqref{PID:SD}, \eqref{ETCon}, \eqref{ETC:PID} with $k_p$, $k_i$, $k_d$ given by \eqref{gains} and \change{$q=\lfloor h^{-\frac12}\rfloor$} exponentially stabilizes the system \eqref{PIDplant1} with any given decay rate $\alpha<\alpha'$ if the sampling period $h>0$ and the event-triggering threshold $\sigma\in[0,1)$ are small enough. \label{th:3:ii}
	\end{itemize}
\end{theorem}

{\em Proof} is given in Appendix~\ref{proof:th:3}. 

\begin{remark}\label{rem:sigma=0}
The event-triggered control \eqref{PID:SD}, \eqref{ETCon}, \eqref{ETC:PID} with $\sigma=0$ degenerates into sampled-data control~\eqref{PID:SD}. Therefore, Theorem~\ref{th:3} with $\sigma=0$ gives the stability conditions for the system \eqref{PIDplant} under the sampled-data PID control \eqref{PID:SD}. 
\end{remark}
\begin{remark}
	Appropriate values of $h$ and $\sigma$ can be found in a manner similar to Remarks~\ref{rem:h} and \ref{rem:sigma}. 
\end{remark}
\subsection{Example}
Following \cite{Ramirez2016}, we consider \eqref{PIDplant} with $a_1=8.4$, $a_2=0$, $b=35.71$. The system is not asymptotically stable if $u=0$. The PID controller \eqref{PID} with $\bar{k}_p=-10$, $\bar{k}_i=-40$, $\bar{k}_d=-0.65$ exponentially stabilizes it with the decay rate $\alpha'\approx10.4$. 

Theorem~\ref{th:3} with $\sigma=0$ (see Remark~\ref{rem:sigma=0}) guarantees that the sampled-data PID controller \eqref{PID:SD} can achieve any decay rate $\alpha<\alpha'$ if the sampling period $h>0$ is small enough. Since $\alpha'$ is on the verge of stability, $\alpha$ close to $\alpha'$ requires to use small $h$. Thus, for $\alpha=10.3$, the LMIs of Theorem~\ref{th:3} are feasible with $h\approx 10^{-7}$, $q=4272$, and $k_p$, $k_i$, $k_d$ given by \eqref{PIDk}. To avoid small sampling \change{period}, we take $\alpha=5$. 

For $\sigma=0$, $\alpha=5$ and each $q=1,2,3,\ldots$ we find the maximum sampling period $h>0$ such that the LMIs of Theorem~\ref{th:3} are feasible. The largest $h$ corresponds to 
\begin{equation*}
\begin{array}{c}
\alpha=5,\quad\sigma=0,\quad q=7,\quad h=4.7\times 10^{-3},\\
k_p\approx-29.76,\quad k_i=-40,\quad k_d\approx19.76, 
\end{array}
\end{equation*}
where $k_p$, $k_i$, $k_d$ are calculated using \eqref{PIDk}. \change{Remark~\ref{rem:sigma=0} implies that} the sampled-data PID controller \eqref{PID:SD} stabilizes \eqref{PIDplant}. 

Theorem~\ref{th:3} remains feasible for 
\begin{equation*}
\begin{array}{c}
\alpha=5,\quad\sigma=9\times 10^{-3},\quad q=7,\quad h=4\times 10^{-3},\\
k_p\approx-33.21,\quad k_i=-40,\quad k_d\approx23.21, 
\end{array}
\end{equation*}
where $k_p$, $k_i$, $k_d$ are calculated using \eqref{PIDk}. Thus, the event-triggered PID control \eqref{PID:SD}, \eqref{ETCon}, \eqref{ETC:PID} exponentially stabilizes~\eqref{PIDplant1}. Performing numerical simulations in a manner described in Section~\ref{sec:tripint}, we find that the event-triggered PID control requires to transmit on average $628.4$ control signals during $10$ seconds. The sampled-data controller \eqref{PID:SD} requires $\left\lfloor\frac{10}{h}\right\rfloor+1=2128$ transmissions. Thus, the event-triggering mechanism reduces the workload of the controller-to-actuator network by more than $70\%$. The total workload of both networks is reduced by more than $26	\%$.

\appendices
\section{Proof of Theorem~\ref{th:1}}\label{proof:th:1}
(i) Consider the functional 
\begin{equation}\label{V}
V=V_0+V_{\delta0}+V_\delta+V_\kappa
\end{equation}
with
\begin{equation*}\arraycolsep=1.4pt
\begin{array}{rl}
&V_0=x^TPx,\\
&V_{\delta0}=h^2e^{2\alpha h}\int_{t_k}^te^{-2\alpha(t-s)}\dot y^T(s)K_0^TW_0K_0\dot y(s)\,ds\\
&\phantom{V_{\delta0}=}-\frac{\pi^2}{4}\int_{t_k}^{t}e^{-2\alpha(t-s)}[y(s)-y(t_k)]^TK_0^TW_0\times\\
\multicolumn{2}{r}{K_0[y(s)-y(t_k)]\,ds,}\\
&V_\delta=h^2e^{2\alpha h}\sum_{i=1}^{r-1}\int_{t_k-q_ih}^te^{-2\alpha(t-s)}\dot y^T(s)K_i^TW_iK_i\dot y(s)\,ds\\
&\phantom{V_\delta=}-\frac{\pi^2}{4}\sum_{i=1}^{r-1}\int_{t_k-q_ih}^{t-q_ih}e^{-2\alpha(t-s)}[y(s)-y(t_k-q_ih)]^T\times\\
\multicolumn{2}{r}{K_i^TW_iK_i[y(s)-y(t_k-q_ih)]\,ds,}\\
&V_\kappa=\sum_{i=1}^{r-1}\int_{t-q_ih}^te^{-2\alpha(t-s)}(s-t+q_ih)^r\times\\
\multicolumn{2}{r}{(y^{(r)}(s))^TK_i^TR_iK_iy^{(r)}(s)\,ds.}
\end{array}
\end{equation*}
The term $V_\kappa$, introduced in \cite{Fridman2017}, compensates Taylor's remainders $\kappa_i$, while $V_{\delta0}$ and $V_\delta$, introduced in \cite{Liu2012}, compensate the sampling errors $\delta_0$ and $\delta$. The Wirtinger inequality (Lemma~\ref{lem:Wirt}) implies $V_{\delta0}\ge0$ and $V_\delta\ge0$. Using \eqref{ClosedLoop} and \eqref{y'}, we obtain 
\begin{equation*}
\arraycolsep=1.4pt
\begin{array}{l}
\dot V_0+2\alpha V_0\!=\!2x^TP[Dx\!+\!BK_0\delta_0\!+\!BK\delta\!+\!BK\kappa]\!+\!2\alpha x^TPx,\\
\dot V_{\delta0}+2\alpha V_{\delta0}=h^2e^{2\alpha h}x^T(K_0CA)^TW_0(K_0CA)x\\
\multicolumn{1}{r}{-\frac{\pi^2}{4}\delta_0^TK_0^TW_0K_0\delta_0,}\\
\dot V_\delta+2\alpha V_\delta=h^2e^{2\alpha h}\sum_{i=1}^{r-1}x^T(K_iCA)^TW_i(K_iCA)x\\
\multicolumn{1}{r}{-\frac{\pi^2}{4}\sum_{i=1}^{r-1}e^{-2\alpha q_ih}\delta_i^TK_i^TW_iK_i\delta_i.}
\end{array}
\end{equation*}
Using $y^{(r)}=CA^{r-1}\dot x$ (which follows from \eqref{y'}) and Jensen's inequality (Lemma~\ref{lem:Jens}) with $\rho(s)=(s-t+q_ih)^{r-1}$, we have
\begin{equation*}
\arraycolsep=1.4pt
\begin{array}{rl}
\dot V_\kappa+2\alpha V_\kappa=&\sum_{i=1}^{r-1}(q_ih)^r(y^{(r)}(t))^TK_i^TR_iK_iy^{(r)}(t)\\ \textstyle
&-\sum_{i=1}^{r-1}r\int_{t-q_ih}^te^{-2\alpha(t-s)}(s-t+q_ih)^{r-1}\times\\
&\multicolumn{1}{r}{(y^{(r)}(s))^TK_i^TR_iK_iy^{(r)}(s)\,ds}\\
\le&\sum_{i=1}^{r-1}(q_ih)^{r}\dot x^T(A^{r-1})^TC^TK_i^TR_iK_iCA^{r-1}\dot x\\
&\multicolumn{1}{r}{-\sum_{i=1}^{r-1}\frac{(r!)^2}{(q_ih)^{r}}e^{-2\alpha q_ih}\kappa_i^TK_i^TR_iK_i\kappa_i.}
\end{array}
\end{equation*}
Summing up, we obtain
\begin{equation}\label{dotV}
\dot V+2\alpha V\le\varphi^T\bar{\Phi}\varphi+\dot x^T(CA^{r-1})^TH(CA^{r-1})\dot x, 
\end{equation}
where 
\begin{equation}\label{varphi}
\varphi=\operatorname{col}\{x,K_0\delta_0,\ldots,K_{r-1}\delta_{r-1},K_1\kappa_1,\ldots,K_{r-1}\kappa_{r-1}\}
\end{equation}
and $\bar{\Phi}$ is obtained from $\Phi$ by removing the last block-column and block-row. Substituting \eqref{ClosedLoop}
for $\dot x$ and applying the Schur complement, we find that $\Phi\le0$ guarantees $\dot V\le-2\alpha V$. \change{Since $V(t_k)\le V(t_k^-)$, the latter} implies exponential stability of the system \eqref{ClosedLoop} and, therefore, \eqref{LTI}, \eqref{control:discrete}. 

(ii) \change{Since $q_i=O(h^{\frac{1}{r}-1})$}, \cite[Lemma~2.1]{Fridman2017} guaranties $M^{-1}=O(h^{\frac1r-1})$, which implies $K_i=O(h^{\frac1r-1})$ for $i=0,\ldots,r-1$. Since $D=\bar{D}$ (with $\bar{D}$ defined in \eqref{ClosedLoop:derivative}) and \eqref{LTI}, \eqref{control:derivative} is exponentially stable with the decay rate $\alpha'$, there exists $P>0$ such that $PD+D^TP+2\alpha P<0$ for any $\alpha<\alpha'$. Choose $W_0=O(h^{-\frac1r})$, $W_i=O(h^{-\frac1r})$, and $R_i=O(h^{1-\frac1r})$ for $i=1,\ldots,r-1$. Applying the Schur complement to $\Phi\le0$, we obtain 
\begin{equation*}
PD+D^TP+2\alpha P+O(h^{\frac1r})<0,
\end{equation*}
which holds for small $h>0$. Thus, (i) guarantees (ii). 
\section{Proof of Theorem~\ref{th:2}}\label{proof:th:2}
Denote $e_k=\hat u_k-u(t_k)$, $k\in\N_0$. The event-triggering mechanism \eqref{ETC}, \eqref{ETControl} guarantees 
\begin{equation}\label{Sproc}
0\le\sigma u^T(t_k)\Omega u(t_k)-e_k^T\Omega e_k. 
\end{equation}
Substituting $\hat u_k=u(t_k)+e_k$ into \eqref{LTI:ETC} and using \eqref{control:errors}, we obtain (cf.~\eqref{ClosedLoop})
\begin{equation}\label{dotXETC}
\dot x=Dx+BK_0\delta_0+BK\delta+BK\kappa+Be_k,\quad t\in[t_k,t_{k+1})
\end{equation}
with $D$ given in \eqref{ClosedLoop}. 
Consider $V$ from \eqref{V}. Calculations similar to those from the proof of Theorem~\ref{th:1} lead to (cf.~\eqref{dotV})
\begin{equation*}\arraycolsep=1.4pt
\begin{array}{l}
\dot V+2\alpha V\stackrel{\eqref{Sproc}}{\le}\dot V+2\alpha V+\sigma u^T(t_k)\Omega u(t_k)-e_k^T\Omega e_k\\
\le\varphi_e^T\bar\Phi_e\varphi_e+\dot x^T(CA^{r-1})^TH(CA^{r-1})\dot x+\sigma u^T(t_k)\Omega u(t_k), 
\end{array}
\end{equation*}
where $\varphi_e=\operatorname{col}\{\varphi,e_k\}$ (with $\varphi$ from \eqref{varphi}) and $\bar{\Phi}_e$ is obtained from $\Phi_e$ by removing the blocks $\Phi_{ij}$ with $i\in\{4,6\}$ or $j\in\{4,6\}$. Substituting \eqref{dotXETC} for $\dot x$ and \eqref{control:errors} for $u(t_k)$ and applying the Schur complement, we find that $\Phi_e\le0$ guarantees $\dot V\le-2\alpha V$. \change{Since $V(t_k)\le V(t_k^-)$, the latter} implies exponential stability of the system \eqref{dotXETC} and, therefore, \eqref{LTI:ETC} under the controller \eqref{control:discrete}, \eqref{ETC}, \eqref{ETControl}. 
\section{Proof of Theorem~\ref{th:3}}\label{proof:th:3}
(i) Consider the functional 
\begin{equation*}
V=V_0+V_v+V_\delta+V_\kappa 
\end{equation*}
with 
\begin{equation*}
\arraycolsep=1.4pt
\begin{array}{rl}
V_0=&x^TPx,\\
V_v=&h^2e^{2\alpha h}\int_{t_k}^te^{-2\alpha(t-s)}\dot x^T(s)S\dot x(s)\,ds\\
\multicolumn{2}{r}{-\frac{\pi^2}{4}\int_{t_k}^te^{-2\alpha(t-s)}v^T(s)Sv(s)\,ds,}\\
V_\delta=&Wk_d^2h^2e^{2\alpha h}\int_{t_k-qh}^te^{-2\alpha(t-s)}(\dot y(s))^2\,ds\\
&-Wk_d^2\frac{\pi^2}{4}\int_{t_k-qh}^{t-qh}e^{-2\alpha(t-s)}[y(s)-y(t_k-qh)]^2ds,\\
V_\kappa=&R k_d^2\int_{t-qh}^te^{-2\alpha(t-s)}(s-t+qh)^2(\ddot y(s))^2\,ds. \\
\end{array}
\end{equation*}
The Wirtinger inequality (Lemma~\ref{lem:Wirt}) implies $V_v\ge0$ and $V_\delta\ge0$. Using the representation \eqref{plant}, we obtain 
\begin{equation*}
\arraycolsep=1.4pt
\begin{array}{l}
\dot V_0+2\alpha V_0=2x^TP[Ax(t)+A_vv(t)+Bk_d(\kappa+\delta)+Be_k]\\
\multicolumn{1}{r}{+2\alpha x^TPx,}\\
\dot V_v+2\alpha V_v=h^2e^{2\alpha h}\dot x^T(t)S\dot x(t)-\frac{\pi^2}{4}v^T(t)Sv(t),\\
\dot V_\delta+2\alpha V_\delta=Wk_d^2h^2e^{2\alpha h}(x_2(t))^2-Wk_d^2\frac{\pi^2}{4}e^{-2\alpha qh}\delta^2(t). 
\end{array}
\end{equation*}
Using Jensen's inequality (Lemma~\ref{lem:Jens}) with $\rho(s)=(s-t+q_ih)$, we obtain 
\begin{equation*}\arraycolsep=1.4pt
\begin{array}{rl}
\dot V_\kappa+2\alpha V_\kappa=&R k_d^2(qh)^2(\ddot y(t))^2\\
&-R k_d^22\int_{t-qh}^te^{-2\alpha(t-s)}(s-t+qh)(\ddot y(s))^2\,ds\\
\le&R k_d^2(qh)^2(\dot x_2(t))^2-R k_d^2\frac{4}{(qh)^2}e^{-2\alpha qh}\kappa^2(t). 
\end{array}
\end{equation*}
For $\omega\ge0$, the event-triggering rule \eqref{ETCon}, \eqref{ETC:PID} guarantees 
$$0\le \omega\sigma u^2(t_k)-\omega e_k^2.$$ 
Thus, we have 
\begin{multline*}
\textstyle \dot V+2\alpha V\le\dot V+2\alpha V+[\omega\sigma u^2(t_k)-\omega e_k^2]\\
\textstyle \le\psi^T\bar{\Psi}\psi+\dot x^T(t)G\dot x(t)+\omega\sigma u^2(t_k), 
\end{multline*}
where $\psi=\operatorname{col}\{x,v/\sqrt{h},k_d\delta,k_d\kappa,e_k\}$ and $\bar{\Psi}$ is obtained from $\Psi$ by removing the last two block-columns and block-rows. Substituting \eqref{plant} for $\dot x$ and \eqref{uk} for $u(t_k)$ and applying the Schur complement, we find that $\Psi\le0$ guarantees $\dot V\le-2\alpha V$. \change{Since $V(t_k)\le V(t_k^-)$, the latter} implies exponential stability of the system \eqref{plant} and, therefore, \eqref{PID:SD}, \eqref{PIDplant1}--\eqref{ETC:PID}. 

(ii) The closed-loop system \eqref{PIDplant}, \eqref{PID} is equivalent to $\dot x=\bar{A}x$ with 
\begin{equation*}
\bar{A}=\begin{bmatrix}
0 & 1 & 0 \\
-a_2+b\bar{k}_p & -a_1+b\bar{k}_d & b\bar{k}_i \\
1 & 0 & 0
\end{bmatrix},\quad x=\begin{bmatrix}
y \\ \dot y \\ \int_0^ty(s)\,ds
\end{bmatrix}. 
\end{equation*}
\change{Since $q=O(\frac{1}{\sqrt{h}})$}, relations \eqref{PIDk} imply $k_p=O(\frac{1}{\sqrt{h}})$, $k_d=O(\frac{1}{\sqrt{h}})$. Since \eqref{PIDplant}, \eqref{PID} is exponentially stable with the decay rate $\alpha'$ and \eqref{PIDk} implies $A=\bar{A}$, there exists $P>0$ such that $PA+A^TP+2\alpha P<0$ for any $\alpha<\alpha'$. Choose $S=O(\frac{1}{h\sqrt{h}})$, $W=O(\frac1{\sqrt{h}})$, $R=O(\sqrt{h})$, and $\omega=O(\frac1{\sqrt{h}})$. Applying the Schur complement to $\Psi\le0$, we obtain 
\begin{equation*}
PA+A^TP+2\alpha P+O(\sqrt{h})+\sigma F<0 
\end{equation*}
with some $F$ independent of $\sigma$. The latter holds for small $h>0$ and $\sigma\ge0$. Thus, (i) guarantees (ii). 
\ifCLASSOPTIONcaptionsoff
  \newpage
\fi

\bibliographystyle{IEEEtran}
\bibliography{library}

\end{document}